\documentclass[a4paper,11pt]{article}

\usepackage{fullpage,amsthm}
\newtheoremstyle{break}
  {0em}{0em}%
  {\itshape}{}%
  {\bfseries}{}%
  {\newline}{}%
\theoremstyle{break}
\usepackage[linewidth=0.5pt,backgroundcolor=frame,linecolor=white]{mdframed}
\newmdtheoremenv[linewidth=0.5pt,backgroundcolor=frame,linecolor=white,innertopmargin=0.4em,innerbottommargin=0.4em,innerleftmargin=0.4em,innerrightmargin=0.4em]{definition}{Definition}
\newmdtheoremenv[linewidth=0.5pt,backgroundcolor=frame,linecolor=white,innertopmargin=0.4em,innerbottommargin=0.4em,innerleftmargin=0.4em,innerrightmargin=0.4em]{algorithm}{Algorithm}
\definecolor{frame}{rgb}{0.9,0.9,0.9}
\definecolor{colorgray}{rgb}{0.8,0.8,0.8}
\usepackage{enumerate}
\usepackage{amsmath,subfig,float,graphicx,multirow}
\usepackage{hyperref}
\hypersetup{colorlinks = true, linktocpage = true, bookmarksopen = true, linkcolor = blue, urlcolor=blue, citecolor = blue, allcolors=blue}
\graphicspath{{"./Figures/"}}
\usepackage[numbers,sort&compress]{natbib}

\begin{document}

\title{Calculating how long it takes for a diffusion process to effectively reach steady state without computing the transient solution}
%\date{Version 1: 21 April 2017\\ Version 2: \today}
\date{\small Submitted on 21 April 2017, Revised 11 July 2017}
\author{Elliot~J.~Carr\\ \small School of Mathematical Sciences, Queensland University of Technology (QUT), Brisbane, Australia\\\small \href{mailto:elliot.carr@qut.edu.au}{elliot.carr@qut.edu.au}}%\\ Brisbane, Australia\\ elliot.carr@qut.edu.au}
\maketitle

\begin{abstract}
Mathematically, it takes an infinite amount of time for the transient solution of a diffusion equation to transition from initial to steady state. Calculating a \textit{finite} transition time, defined as the time required for the transient solution to transition to within a small prescribed tolerance of the steady state solution, is much more useful in practice. In this paper, we study estimates of finite transition times that avoid explicit calculation of the transient solution by using the property that the transition to steady state defines a cumulative distribution function when time is treated as a random variable. In total, three approaches are studied: (i) mean action time (ii) mean plus one standard deviation of action time and (iii) a new approach derived by approximating the large time asymptotic behaviour of the cumulative distribution function. The new approach leads to a simple formula for calculating the finite transition time that depends on the prescribed tolerance $\delta$ and the $(k-1)$th and $k$th moments ($k \geq 1$) of the distribution. Results comparing exact and approximate finite transition times lead to two key findings. Firstly, while the first two approaches are useful at characterising the time scale of the transition, they do not provide accurate estimates for diffusion processes. Secondly, the new approach allows one to calculate finite transition times accurate to effectively any number of significant digits, using only the moments, with the accuracy increasing as the index $k$ is increased. 
\end{abstract}

\section{Introduction}
A common question of practical interest that arises when modelling transport phenomena is \textit{how long does the process take}? For example, how long does it take for a bar to heat up or how long does it take for an initial concentration field to distribute evenly? Mathematically, the correct answer is the highly impractical answer: it takes an infinite amount of time for the transient solution to transition from initial to steady state. 

Finite (and hence practically useful) answers to the above question have been studied and proposed by many authors \cite{hickson_thesis_2010,hickson_2011a,hickson_2009_part1,hickson_2009_part2,simpson_2013,ellery_2012b,landman_2000,mcnabb_1991,mcnabb_1993}. In this paper, we consider a finite transition time defined as the amount of time required for the transient solution to transition from its initial state to within a prescribed tolerance of its steady state \cite{simpson_2013}. To illustrate this definition, consider a transport process on an interval $\mathcal{L} := (l_{0},l_{m})$ and let $u(x,t)$ be the transient solution at position $x\in\mathcal{L}$ and time $t> 0$, $u_{0}(x)$ be the given initial solution and $u_{\infty}(x)$ be the steady state solution. We make the following assumptions:
\begin{enumerate}[(i)]
\item $u_{0}(x)\neq u_{\infty}(x)$ for all $x\in\mathcal{L}$, 
\item For a given $x\in\mathcal{L}$, $u(x,t)$ monotonically increases to $u_{\infty}(x)$ (implying $u(x,t) < u_{\infty}(x)$) for all $t > 0$ if $u_{0}(x) < u_{\infty}(x)$ or $u(x,t)$ monotonically decreases to $u_{\infty}(x)$ (implying $u(x,t) > u_{\infty}(x)$) for all $t > 0$ if $u_{0}(x)>u_{\infty}(x)$,
\item $u_{0}(x)$ is piecewise continuous on $\mathcal{L}$.
\end{enumerate}

In this work, we distinguish between \textit{local} and \textit{global} transition times as follows. The local transition time provides a finite measure of the time required to transition from initial to steady state at position $x$. On the other hand, the global transition time is a finite measure of the time required for the entire transport process to effectively reach steady state, and hence provides a practical finite answer to the question of how long the process takes. Both of these quantities are defined below, where we remark that the assumption (ii) listed above ensures that the left-hand side of equation (\ref{eq:singlelayer_local_transition_time}) is always positive. 

\begin{definition}[Local transition time]
\label{def:singlelayer_local_transition_time}
\normalfont
The~\textit{local transition time},~denoted by $t_{s}(x)$, provides the transition time as a function of position $x\in\mathcal{R} := \{x\in\mathcal{L}\,|\,u_{0}(x)\neq u_{\infty}(x)\}$, and is defined as the value of $t > 0$ satisfying:
\begin{equation}
\label{eq:singlelayer_local_transition_time}
\frac{u(x,t)-u_{\infty}(x)}{u_{0}(x)-u_{\infty}(x)} = \delta,
\end{equation}
where $0<\delta\ll 1$ is a specified tolerance.% (i.e., $\delta = 10^{-p}$, $p=1,2,\hdots$). 
\end{definition}%\vspace*{-0.3cm}

\begin{definition}[Global transition time]
\label{def:singlelayer_global_transition_time}
\normalfont
The \textit{global transition time} is defined as the supremum of the local finite transition time:
\begin{equation}
\label{eq:singlelayer_global_transition_time}
\widehat{t}_{s} := \sup_{x\in \mathcal{R}}\,t_{s}(x)
\end{equation}
where the supremum is taken over the domain of $t_{s}(x)$.
\end{definition}%\vspace*{-0.3cm}

%The goal of this paper is to estimate the above local and global transition times without explicit calculation of the transient solution $u(x,t)$. 

A common finite measure of the time required to reach steady state is the \textit{mean action time} \cite{simpson_2013,jazaei_2014} introduced by \citet{mcnabb_1991}. The mean action time takes a probabilistic approach to the problem by treating time $t$ as continuous random variable with support $[0,\infty)$ and computing its mean or expected value. This is achieved by utilising the observation that the function
\begin{equation}
\label{eq:singlelayer_CDF}
F(t; x) := 1 - \frac{u(x,t)-u_{\infty}(x)}{u_{0}(x)-u_{\infty}(x)},
\end{equation}
satisfies $F(0; x) = 0$ and $\lim_{t\rightarrow\infty}F(t; x) = 1$ and hence defines a cumulative distribution function of $t$, parametrised in terms of the position $x$. Using the cumulative distribution function, the local transition time (Definition \ref{eq:singlelayer_local_transition_time}) can now be reformulated as follows:
\begin{gather}
\label{eq:transition_time_x_Ftx}
\textit{$t_{s}(x)$ is the value of $t > 0$ satisfying $F(t; x) = 1- \delta$}.
\end{gather}
Note that equation (\ref{eq:transition_time_x_Ftx}) is now equivalent to a classical inverse problem in probability theory: find $t_{s}(x)$ such that $P(t\geq t_{s}(x)) = \delta$. 

Differentiating (\ref{eq:singlelayer_CDF}) with respect to $t$ yields the probability density function \cite{ellery_2012a,simpson_2013,simpson_2017}:
\begin{equation}
\label{eq:singlelayer_PDF}
f(t; x) = \frac{1}{u_{\infty}(x)-u_{0}(x)}\frac{\partial}{\partial t}\left[u(x,t) - u_{\infty}(x)\right],
\end{equation}
which allows the $k$th moment, at position $x$, to be defined as follows:
\begin{equation}
\label{eq:singlelayer_raw_moment}
M_{k}(x) := \int_{0}^{\infty} t^{k}f(t;x)\,dt.
\end{equation}
The mean action time at position $x$, denoted henceforth by $\text{MAT}(x)$, is then defined by the mean or first moment $M_{1}(x)$. 

The attractiveness of the mean action time is that it is possible to obtain the explicit form of $\text{MAT}(x)$, and hence an estimate of the finite amount of time required for the transient solution to effectively reach steady state at position $x$, without explicitly calculating the transient solution $u(x,t)$ \cite{ellery_2012b,simpson_2013,simpson_2017}. For this reason, the mean action time is a popular transition time estimate that has been used in several applications including freezing and thawing \cite{landman_2000}, morphogen formation in the formation of tissues and organs \cite{berezhkovskii_2010, simpson_2017} and groundwater modelling \cite{simpson_2013,jazaei_2014,jazaei_2017}. Another argument in favour of the mean action time is that it does not suffer from the subjectivity of choosing a threshold tolerance \cite{ellery_2012b}, as in Definition \ref{def:singlelayer_local_transition_time}. While this is true, the question of how close the transient solution, evaluated at the supremum of the mean action time, is to the steady state still remains and involves some discretion. The opinion of this author is that the tolerance $\delta$ is useful and should be inferred by the physical problem under consideration. For example, in many problems involving comparison to experimental data, a suitable choice for $\delta$ is the error associated with the measurement device.

Since McNabb and Wake's \cite{mcnabb_1991}, notable contributions to the theory of mean action time, almost exclusively focussed on one-dimensional linear homogeneous problems, have appeared in several papers  \cite{landman_2000,hickson_2011b,ellery_2012b,ellery_2012a,simpson_2013,jazaei_2014,simpson_2017,jazaei_2017}. For discussion on higher-dimensional, nonlinear and/or heterogeneous problems, which do not form the focus of this work, the reader is directed to three papers: \cite{hickson_2011b,ellery_2012b,landman_2000}.

Research on mean action time was revived in the early 2010s by both \citet{berezhkovskii_2010}, who proposed the definition of \textit{local accumulation time} for a reaction-diffusion model related to the study of morphogen formation and \citet{ellery_2012b}, who demonstrated that this definition was equivalent to McNabb's mean action time. Over two papers, \citet{ellery_2012a,ellery_2012b} derived the mean action time for a linear advection-diffusion-reaction process and presented a framework for computing the higher central moments (called \textit{moments of action} by the authors), demonstrating how in each case exact expressions can be found without explicit calculation of the transient solution. The second moment of action, the \textit{variance of action time}, was noted as being particularly insightful as a small value implies that the mean action time is a useful estimate of the time required to effectively reach steady state \cite{ellery_2012b,simpson_2013}. \citet{simpson_2013} applied and extended the theory to a linearised groundwater flow model, governed by a reaction-diffusion equation with constant source term, to study the time required for a transient response (such as an aquifer recharge/discharge process) to effectively reach steady state. Comparing the model predictions to a laboratory scale experiment, it was found that while the mean action time underestimates the time required to transition to steady state, significant improvement can be obtained by adding one standard deviation of action time. \citet{jazaei_2014} extended the derivations for groundwater flow to time-dependent boundary conditions, where it was again argued, via visual observation of the experimental data, that the \textit{mean plus one standard deviation of action time} was a good approximation for the time required for the system to effectively reach steady-state. Most recently, \citet{simpson_2017} and \citet{jazaei_2017} have extended the theory to estimate the time required for the gradient or flux of the transient solution to approach their corresponding steady-state values. 

In this paper, we study transition time estimates for a linear homogeneous diffusion equation with general time-independent boundary conditions. Two key contributions to the literature are presented, namely we:
\begin{enumerate}[(i)]
\item demonstrate that, while both the mean action time and the mean plus one standard deviation of action time are useful at characterising the associated time-scale, they do not provide accurate estimates of the transition time for diffusion processes. 
\item propose a highly accurate alternative approach for estimating transition times based on using higher-order consecutive moments to approximate the large time asymptotic behaviour of the cumulative distribution function $F(t; x)$. 
\end{enumerate}
The second item is the main contribution of the work, with the 
novel approach providing a simple highly-accurate formula for calculating the local transition time (\ref{eq:singlelayer_local_transition_time}) (and hence global transition time (\ref{eq:singlelayer_global_transition_time})) using only two consecutive moments $M_{k-1}(x)$ and $M_{k}(x)$ ($k \geq 1$). The result links, for the first time, the higher-order moments to calculation of finite transition times. This represents a significant breakthrough as previously it was believed that such calculation required the transient solution $u(x,t)$ \cite{ellery_2012b,simpson_2013}.

\section{Diffusion model}
\label{sec:singlelayer_diffusion}
The diffusion model considered in this paper is described below. Consider a diffusion process on the interval $\mathcal{L} := [l_{0},l_{m}]$ governed by the linear homogeneous diffusion equation:
\begin{subequations}
\label{eq:singlelayer_problem}
\begin{gather}
\label{eq:singlelayer_transient_PDE}
\frac{\partial u}{\partial t} = D\frac{\partial^{2}u}{\partial x^{2}},
\end{gather}
for $x\in\mathcal{L}$ and $t > 0$, subject to the initial condition:
\begin{gather}
\label{eq:singlelayer_transient_IC}
u(x,0) = u_{0}(x),
\end{gather}
for $x\in\mathcal{L}$ and boundary conditions:
\begin{gather}
\label{eq:singlelayer_transient_BC1}
a_{L}u(l_{0},t) - b_{L}\frac{\partial u}{\partial x}(l_{0},t) = c_{L},\\
\label{eq:singlelayer_transient_BC2}
a_{R}u(l_{m},t) + b_{R}\frac{\partial u}{\partial x}(l_{m},t) = c_{R},
\end{gather}
for $t > 0$. 
\end{subequations}
In the equations listed above, $u(x,t)$ is the solution (e.g. temperature/concentration) at position $x$ and time $t$, $u_{0}(x)$ is the specified initial condition and $D>0$ is the constant diffusion coefficient. In the boundary conditions, $a_{L}$, $b_{L}$, $c_{L}$, $a_{R}$, $b_{R}$ and $c_{R}$ are constants satisfying: $a_{L}\geq 0$, $b_{L}\geq 0$, $a_{R}\geq 0$, $b_{R}\geq 0$, $a_{L} + b_{L} > 0$ and $a_{R} + b_{R} > 0$ \cite{carr_2017a,carr_2017b}.

The corresponding steady-state solution of (\ref{eq:singlelayer_problem}), denoted by $u_{\infty}(x)$, is the linear function satisfying the boundary value problem:
\begin{subequations}
\label{eq:singelayer_steady_state_problem}
\begin{gather}
\label{eq:singelayer_steady_state_DE}
Du_{\infty}''(x) = 0,\quad x\in \mathcal{L},\\
\label{eq:singlelayer_steady_state_BC1}
a_{L}u_{\infty}(l_{0}) - b_{L}u_{\infty}'(l_{0}) = c_{L},\\
\label{eq:singlelayer_steady_state_BC2}
a_{R}u_{\infty}(l_{m}) + b_{R}u_{\infty}'(l_{m}) = c_{R}.
\end{gather}
\end{subequations}
For the Neumann problem ($u_{\infty}'(l_{0})=u_{\infty}'(l_{m})=0$), it is well-known \cite{carr_2016a} that the solution of (\ref{eq:singelayer_steady_state_problem}) is unique only up to an additive constant and that including the additional constraint:
\begin{gather}
\label{eq:singlelayer_steady_state_constraint}
\int_{l_{0}}^{l_{m}} u_{\infty}(x)\,dx = \int_{l_{0}}^{l_{m}} u_{0}(x)\,dx,
\end{gather}
gives the correct steady state solution of (\ref{eq:singlelayer_problem}). %In the analysis that follows, $u_{\infty}(x)$ is assumed known.

\section{Computing the moments}
In this section, we present an algorithm for calculating the first $q$ moments (\ref{eq:singlelayer_raw_moment}), namely $M_{k}(x)$ ($k = 1,\hdots,q$), for the diffusion problem (\ref{eq:singlelayer_problem}). The derivation follows closely the procedure taken by several authors \cite{ellery_2012a,ellery_2013b,simpson_2017} with the exception that we are interested in the raw moments as opposed to the central moments.

Combining equations (\ref{eq:singlelayer_PDF}) and (\ref{eq:singlelayer_raw_moment}) gives the following expression for the $k$th moment:
\begin{align}
\label{eq:singlayer_M_2a}
M_{k}(x) = \frac{1}{h(x)}\int_{0}^{\infty} t^{k}\frac{\partial}{\partial t}\left[u(x,t)-u_{\infty}(x)\right]\,dt,
\end{align} 
for $k = 0,1,\hdots$,~where $h(x):=u_{\infty}(x)-u_{0}(x)$. Applying integration by parts and noting that $\lim\limits_{t\rightarrow\infty} t^{k}[u(x,t)-u_{\infty}(x)] = 0$ \cite{ellery_2013b} yields:
\begin{gather}
\label{eq:singlelayer_M_2}
M_{k}(x) = \frac{k}{h(x)}\int_{0}^{\infty} t^{k-1}\left[u_{\infty}(x)-u(x,t)\right]\,dt.
\end{gather}
For all integers $k > 1$, each of the above moments can be calculated without requiring the transient solution $u(x,t)$ by deriving a boundary value problem satisfied by $M_{k}(x)$ as follows. Define:
\begin{align}
\label{eq:singlelayer_Mbar}
\overline{M}_{k}(x) := M_{k}(x)h(x) = k\int_{0}^{\infty} t^{k-1}\left[u_{\infty}(x)-u(x,t)\right]\,dt,
\end{align}
and consider the derivatives:
\begin{gather}
\label{eq:singlelayer_Mbar_derivative}
\overline{M}_{k}'(x) = k\int_{0}^{\infty} t^{k-1}\left[u_{\infty}'(x)-\frac{\partial u}{\partial x}(x,t)\right]\,dt,\\
\label{eq:singlelayer_Mbar_derivative2}
\overline{M}_{k}''(x) = k\int_{0}^{\infty} t^{k-1}\left[u_{\infty}''(x)-\frac{\partial^{2} u}{\partial x^{2}}(x,t)\right]\,dt.
\end{gather}
Using equations (\ref{eq:singlelayer_transient_PDE}) and (\ref{eq:singelayer_steady_state_DE}) in equation (\ref{eq:singlelayer_Mbar_derivative2}) yields
\begin{gather}
\label{eq:singlelayer_moment_ODE_almost}
\overline{M}_{k}''(x) = \frac{k}{D}\int_{0}^{\infty} t^{k-1}\frac{\partial}{\partial t}\left[u_{\infty}(x)-u(x,t)\right]\,dt.
\end{gather}
The above expressions lead to the following boundary-value problem for $\overline{M}_{k}(x)$:
\begin{subequations}
\label{eq:singlelayer_moment_BVP}
\begin{gather}
\label{eq:singlelayer_moment_ODE}
\overline{M}_{k}''(x) = -\frac{k}{D}\overline{M}_{k-1}(x),\quad x\in\mathcal{L},\\
\label{eq:singlelayer_moment_BC1}
a_{L}\overline{M}_{k}(l_{0}) - b_{L}\overline{M}_{k}'(l_{0}) = 0,\\
\label{eq:singlelayer_moment_BC2}
a_{R}\overline{M}_{k}(l_{m}) + b_{R}\overline{M}_{k}'(l_{m}) = 0,
\end{gather}
\end{subequations}
where the right-hand side of the differential equation (\ref{eq:singlelayer_moment_ODE}) is identified from (\ref{eq:singlelayer_moment_ODE_almost}) using $\overline{M}_{k-1}(x) = M_{k-1}(x)h(x)$ and equation (\ref{eq:singlayer_M_2a})
and the boundary conditions (\ref{eq:singlelayer_moment_BC1}) and (\ref{eq:singlelayer_moment_BC2}) are derived by utilising equations (\ref{eq:singlelayer_Mbar}) and (\ref{eq:singlelayer_Mbar_derivative}) and the boundary conditions for $u(x,t)$ and $u_{\infty}(x)$: (\ref{eq:singlelayer_transient_BC1}), (\ref{eq:singlelayer_transient_BC2}), (\ref{eq:singlelayer_steady_state_BC1}) and (\ref{eq:singlelayer_steady_state_BC2}).% intuitive sense

In this paper, we consider only relatively simple initial condition functions $u_{0}(x)$, where the boundary value problem (\ref{eq:singlelayer_moment_BVP}) can be solved analytically by integration. For more complicated initial conditions, a numerical method could easily be applied to solve (\ref{eq:singlelayer_moment_BVP}). Integrating (\ref{eq:singlelayer_moment_ODE}) yields the general solution:
\begin{align}
\label{eq:singlelayer_moment_general_solution}
\overline{M}_{k}(x) &= G_{k}(x) + c_{k,1} + c_{k,2}x,
\end{align}
where 
\begin{gather}
\label{eq:singlelayer_Gk_funcs}
G_{k}(x) := -\frac{k}{D}\int\int\overline{M}_{k-1}(x)\,dx\,dx.
\end{gather}
The constants $c_{k,1}$ and $c_{k,2}$ are determined by substituting the form of (\ref{eq:singlelayer_moment_general_solution}) into the boundary conditions (\ref{eq:singlelayer_moment_BC1}) and (\ref{eq:singlelayer_moment_BC2}) and solving the resulting linear system:
\begin{align}
\label{eq:singlelayer_linear_system}
\mathbf{A}\mathbf{c} = \mathbf{b},
\end{align}
where $\mathbf{c} = \left[c_{k,1}, c_{k,2}\right]^{T}$ and
\begin{align}
\label{eq:singlelayer_linear_system_matrix}
\mathbf{A} &= \left[\begin{matrix} a_{L} & a_{L}l_{0}-b_{L}\\ a_{R} & a_{R}l_{m}+b_{R} \end{matrix}\right]\\
\label{eq:singlelayer_linear_system_rhs}
\mathbf{b} &= \left[\begin{matrix} b_{L}G_{k}'(l_{0}) - a_{L}G_{k}(l_{0})\\ -b_{R}G_{k}'(l_{m})-a_{R}G_{k}(l_{m})\end{matrix}\right],
\end{align}
with: 
\begin{gather*}
G_{k}'(x) = -\frac{k}{D}\int\overline{M}_{k-1}(x)\,dx.
\end{gather*}
We remark that the matrix $\mathbf{A}$ is the same as the one that appears when using a similar strategy to solve the boundary value problem (\ref{eq:singelayer_steady_state_problem}) for the steady-state solution. For the Neumann problem ($\overline{M}_{k}'(l_{0}) = \overline{M}_{k}'(l_{m}) = 0$), where the solution of the linear system (\ref{eq:singlelayer_linear_system}) (and hence the solution of the boundary value problem (\ref{eq:singlelayer_moment_BVP})) is unique only up to an additive constant, we require an additional constraint on $\overline{M}_{k}(x)$ similar to equation (\ref{eq:singlelayer_steady_state_constraint}). Note that equation (\ref{eq:singlelayer_steady_state_constraint}) together with Neumann boundary conditions implies conservation, that is:
\begin{gather}
\label{eq:singlelayer_conservation}
\int_{l_{0}}^{l_{m}} u(x,t)\,dx = \int_{l_{0}}^{l_{m}} u_{\infty}(x)\,dx,
\end{gather}
for all $t \geq 0$. Integrating (\ref{eq:singlelayer_Mbar}) from $x = l_{0}$ to $x = l_{m}$, reversing the order of integration in the resulting double integral and using (\ref{eq:singlelayer_conservation}) gives the required constraint on the solution of (\ref{eq:singlelayer_moment_BVP}) for the Neumann problem:
\begin{gather*}
\int_{l_{0}}^{l_{m}}\overline{M}_{k}(x)\,dx = 0.
\end{gather*}
Incorporating this constraint into the linear system (\ref{eq:singlelayer_linear_system}) yields the slightly modified form:
\begin{align}
\label{eq:singlelayer_linear_system_matrix_neumann}
\mathbf{A} &= \left[\begin{matrix} 0 & -1\\ 0 & -1\\ l_{m}-l_{0} & \frac{1}{2}(l_{m}^2-l_{0}^2) \end{matrix}\right]\\
\label{eq:singlelayer_linear_system_rhs_neumann}
\mathbf{b} &= \left[\begin{matrix} G_{k}'(l_{0})\\ G_{k}'(l_{m})\\ K_{k}(l_0) - K_{k}(l_m)\end{matrix}\right],
\end{align}
where 
\begin{align*}
K_{k}(x) = \int G_{k}(x)\,dx,
\end{align*}
which has a unique solution provided that $G_{k}'(l_{0}) = G_{k}'(l_{m})$. Once $c_{k,1}$ and $c_{k,2}$ are identified and $\overline{M}_{k}(x)$ is determined, the $k$th moment is computed as follows:
\begin{align}
\label{eq:mean_action_time}
M_{k}(x) &:= \frac{1}{h(x)}\left(G_{k}(x) + c_{k,1} + c_{k,2}x\right).
\end{align}
%which 
Noting that $M_{0}(x) = 1$ and hence $\overline{M}_{0}(x) = h(x) = u_{\infty}(x) - u_{0}(x)$, allows the $G_{k}(x)$ functions (\ref{eq:singlelayer_Gk_funcs}) and the moments to be calculated recursively, as outlined in Algorithm \ref{alg:singlelayer_moments}.

%\begin{figure}
%\begin{algorithm}[Single layer moments]
%\label{alg:singlelayer_moments}
%\normalfont
%\begin{flalign*}
%&\text{$\overline{M}_{0}(x) := u_{\infty}(x) - u_{0}(x)$} &\\
%&\text{\textbf{for} $k = 1,\hdots,q$} &\\
%&\text{\qquad $G_{k}'(x) := -\frac{k}{D}\int\overline{M}_{k-1}(x)\,dx$} &\\
%&\text{\qquad $G_{k}(x) := \int G_{k}'(x)\,dx$} &\\
%&\text{\qquad Compute $c_{k,1}$ and $c_{k,2}$ by solving the} &\\
%&\text{\qquad linear system defined by equations (\ref{eq:singlelayer_linear_system}), (\ref{eq:singlelayer_linear_system_matrix})} &\\
%&\text{\qquad and (\ref{eq:singlelayer_linear_system_rhs}) or in the case of Neumann} &\\
%&\text{\qquad boundary conditions the linear system} &\\
%&\text{\qquad defined by equations (\ref{eq:singlelayer_linear_system}), (\ref{eq:singlelayer_linear_system_matrix_neumann}) and (\ref{eq:singlelayer_linear_system_rhs_neumann}).} &\\
%&\text{\qquad $\overline{M}_{k}(x) := G_{k}(x) + c_{k,1}+c_{k,2}x$} &\\
%&\text{\qquad $M_{k}(x) := \frac{\overline{M}_{k}(x)}{\overline{M}_{0}(x)}$} &\\
%&\textbf{end}
%\end{flalign*}
%\end{algorithm}
%\end{figure}
\begin{figure}
\begin{algorithm}[Moments]
\label{alg:singlelayer_moments}
\normalfont
\begin{tabular}{p{0.4cm}p{7.5cm}}
&\\
\multicolumn{2}{l}{$\overline{M}_{0}(x) := u_{\infty}(x) - u_{0}(x)$}\\
\textbf{for} & $k = 1,\hdots,q$\\
 & $G_{k}'(x) := -\frac{k}{D}\int\overline{M}_{k-1}(x)\,dx$\\
 & $G_{k}(x) := \int G_{k}'(x)\,dx$\\
 & Compute $c_{k,1}$ and $c_{k,2}$ by solving the linear system defined by equations (\ref{eq:singlelayer_linear_system}), (\ref{eq:singlelayer_linear_system_matrix}) and (\ref{eq:singlelayer_linear_system_rhs}) or in the case of Neumann boundary conditions the linear system defined by equations (\ref{eq:singlelayer_linear_system}), (\ref{eq:singlelayer_linear_system_matrix_neumann}) and (\ref{eq:singlelayer_linear_system_rhs_neumann}).\\
& $\overline{M}_{k}(x) := G_{k}(x) + c_{k,1}+c_{k,2}x$\\
& $M_{k}(x) := \overline{M}_{k}(x)/\overline{M}_{0}(x)$\\
\textbf{end}
\end{tabular}
\end{algorithm}
\end{figure}

\section{Transition time estimation}

Using the moments, derived in the previous section, we now present three estimates of the local transition time (Definition \ref{def:singlelayer_local_transition_time}), labelled $t_{s}^{(n)}(x)$ for $n = 1,2,3$. For each estimate, the corresponding estimate of the global transition time is defined according to Definition \ref{def:singlelayer_global_transition_time} as:
\begin{gather}
\label{eq:singlelayer_global_transition_time_estimate}
\widehat{t}_{s}^{(n)} := \sup_{x\in\mathcal{R}}t_{s}^{(n)}(x).
\end{gather}

\subsection{Low accuracy estimates using the first and second moments}
\label{sec:singlelayer_low_accuracy_estimates}

\subsubsection{Mean action time}
\label{sec:singlelayer_mean_action_time}
The mean action time $\text{MAT}(x)$, which is often used as an estimate of the time required to reach steady state \cite{mcnabb_1991,landman_2000,simpson_2013}, defines the following estimate of the local transition time (Definition \ref{def:singlelayer_local_transition_time}):
\begin{gather}
\label{eq:singlelayer_local_transition_time_1}
t_{s}^{(1)}(x) := \text{MAT}(x) = M_{1}(x).
\end{gather}
Although $t_{s}^{(1)}(x)$ and the corresponding estimate of the global transition time $\widehat{t}_{s}^{(1)}$ do not depend on the tolerance $\delta$, it will be interesting to investigate how closely the transient solution $u(x,t)$, evaluated at $t = \widehat{t}_{s}^{(1)}$, is to the steady state solution. 

\subsubsection{Mean plus one standard deviation of action time}
\label{sec:singlelayer_mean_plus_standard_deviation}
For most problems, it is probably unreasonable to expect that $\widehat{t}_{s}^{(1)}$ is an accurate estimate of the time required to reach steady state since its unlikely that the cumulative distribution function $F(t; x)$ evaluated at $t = \widehat{t}_{s}^{(1)}(x)$ is close to one\footnote{Or equivalently that \protect{$P(t\geq\widehat{t}_{s}^{(1)}(x))$} is small.} for all $x\in\mathcal{R}$ since $t_{s}^{(1)}(x)$ is the mean of the distribution. \citet{simpson_2013} suggest the mean action time plus one standard deviation of action time (square root of the variance of action time) as a way to improve estimation. Hence, we will also investigate the following estimate of the local transition time (Definition \ref{def:singlelayer_local_transition_time}): 
\begin{align}
\label{eq:singlelayer_local_transition_time_2}
t_{s}^{(2)}(x) := \text{MAT}(x) + \sqrt{\text{VAT}(x)} = M_{1}(x) + \sqrt{M_{2}(x)-M_{1}(x)^2},
\end{align}
which has been equivalently expressed in terms of the first and second moments. Analogously to $t_{s}^{(1)}$, equation (\ref{eq:singlelayer_local_transition_time_2}) also does not depend on the tolerance $\delta$. We remark that the choice of one standard deviation is subjective and there is no reason why any number of standard deviations could not be used. However, we will not pursue this further as, unlike for the normal distribution, in general one cannot determine the probability that a random variable is greater than the sum of its mean and a given number of standard deviations.

\subsection{High accuracy estimates using higher order moments}
\label{sec:singlelayer_high_accuracy_estimates}

To estimate transition times, an accurate approximation of the cumulative distribution function $F(t; x)$ is required only for large $t$ (relative to the problem) where $u(x,t)\approx u_{\infty}(x)$. The exact solution of the diffusion problem (\ref{eq:singlelayer_problem}) has the following functional form \cite{trim_1990,strauss_1992}:
\begin{gather}
\label{eq:singlelayer_solution}
u(x,t) = u_{\infty}(x) + \sum_{n=1}^{\infty}\gamma_{n}(x)e^{-t\xi_{n}},
\end{gather}
where $\gamma_{n}$ and $\xi_{n} > 0$ depend on the eigenvalues and eigenfunctions of the transient solution. Inserting (\ref{eq:singlelayer_solution}) into (\ref{eq:singlelayer_CDF}) it follows that the cumulative distribution function $F(t; x)$ has the functional form:
\begin{gather}
\label{eq:CDF_form}
F(t; x) =  1 - \sum_{n=1}^{\infty}\zeta_{n}(x)e^{-t\xi_{n}},
\end{gather}
where $\zeta_{n}(x) = \gamma_{n}(x) / (u_{0}(x) - u_{\infty}(x))$. Assuming the $\xi$'s are arranged in ascending order (i.e., $\xi_{1} < \xi_{2} < \hdots $) and $\zeta_{1}(x)\neq 0$, then the cumulative distribution function satisfies the following asymptotic relation:
\begin{align}
\label{eq:CDF_large_t}
F(t; x) \sim 1 - \zeta_{1}(x)e^{-t\xi_{1}}\quad\text{for large $t$}.
\end{align}
The above analysis suggests an approximation to the large time behaviour of $F(t; x)$ should be sought in the form given below: 
\begin{align}
\label{eq:CDF_large_t_alpha_beta}
F(t; x) \simeq 1 - \alpha(x) e^{-t\beta(x)}\quad\text{for large $t$},
\end{align}
where $\alpha(x)$ and $\beta(x)$ are as yet unspecified functions. 

If we can devise a method for accurately computing values of $\alpha$ and $\beta$ that are close to $\zeta_{1}$ and $\xi_{1}$ then highly accurate transition time estimates can be calculated. The aim being, of course, to achieve this without explicit calculation of the transient solution $u(x,t)$ (i.e. without explicit computation of any eigenvalues or eigenfunctions) and using only the moments: $M_{k}(x)$ for $k = 1,\hdots,q$. The probability density function corresponding to the cumulative distribution function (\ref{eq:CDF_form}) is given via differentiation:
\begin{align*}
f(t;x) &= \sum_{n=1}^{\infty} \zeta_{n}(x)\xi_{n}e^{-t\xi_{n}}.
\end{align*}
Using this form of $f(t;x)$ in equation (\ref{eq:singlelayer_raw_moment}), the $k$th moment is found to be:
\begin{align*}
M_{k}(x) = k!\sum_{n=1}^{\infty}\frac{\zeta_{n}}{\xi_{n}^{k}}.
\end{align*}
Since $\xi_{1} < \xi_{n}$ for all $n = 2,3,\hdots$ it follows that $M_{k}(x)$ satisfies the asymptotic relation:
\begin{align}
\label{eq:singlelayer_moment_asymptotic}
M_{k}(x) &\sim \frac{\zeta_{1}k!}{\xi_{1}^{k}}\quad\text{for large $k$}.
\end{align}
This latter observation motivates the following pair of coupled equations satisfied by $\alpha(x)$ and $\beta(x)$, formulated by matching the $(k-1)$th and $k$th moments:
\begin{align*}
\frac{\alpha(x)}{\beta(x)^{k-1}} &= \frac{M_{k-1}(x)}{(k-1)!},\\
\frac{\alpha(x)}{\beta(x)^{k}} &= \frac{M_{k}(x)}{k!}.
\end{align*}
Provided both $M_{k-1}(x)\neq 0$ and $M_{k}(x)\neq 0$, the above system of equations can be solved exactly to obtain the following explicit formulae:
\begin{align*}
\alpha_{k}(x) &= \frac{M_{k}(x)}{k!}\left(\frac{kM_{k-1}(x)}{M_{k}(x)}\right)^{k},\\
\beta_{k}(x) &= \frac{kM_{k-1}(x)}{M_{k}(x)},% and $\beta_{1} = e^{\widetilde{\beta}_{1}}$}
\end{align*}
where we have included the subscript $k$ on $\alpha$ and $\beta$ to denote dependence on the $(k-1)$th and $k$th moments. We remark that $M_{k}(x) > 0$ and hence $\alpha_{k}(x) > 0$ and $\beta_{k}(x) > 0$. The analysis above leads to the following approximation of the large time behaviour of the exact cumulative distribution function (\ref{eq:CDF_form}):
\begin{gather}
\label{eq:CDF_moments}
F(t; x) \simeq 1 - \alpha_{k}(x)e^{-t\beta_{k}(x)}\quad\text{for large $t$},
\end{gather}
which involves only the moments $M_{k-1}(x)$ and $M_{k}(x)$. Following definition (\ref{eq:transition_time_x_Ftx}), equating the right-hand side of (\ref{eq:CDF_moments}) with $1 - \delta$ and solving for $t$ yields the following local transition time estimate:
\begin{gather}
\label{eq:singlelayer_local_transition_time_3}
t_{s}^{(3)}(x) := \frac{1}{\beta_{k}(x)}\log\left(\frac{\alpha_{k}(x)}{\delta}\right),
\end{gather}
where $\log$ is the natural (base $e$) logarithm. Note that $t_{s}^{(3)}(x)$ increases with decreasing tolerance $\delta$, which is consistent with the fact that it takes a longer amount of time for the transient solution to transition to within a smaller tolerance of its steady state. Moreover, we require $\delta\geq\alpha_{k}(x)$ for all $x\in \mathcal{R}$ to ensure the obvious physical constraint $t_{s}^{(3)}(x) \geq 0$ is satisfied. Inserting the expressions for $\alpha_{k}(x)$ and $\beta_{k}(x)$ into (\ref{eq:singlelayer_local_transition_time_3}), we obtain the following simple formula for estimating the local transition time depending on the $(k-1)$th and $k$th moments:
\begin{gather}
\label{eq:singlelayer_local_transition_time_3b}
\hspace*{-0.3cm}t_{s}^{(3)}(x) := \frac{M_{k}(x)}{kM_{k-1}(x)}\log\left[\frac{M_{k}(x)}{k! \,\delta}\left(\frac{kM_{k-1}(x)}{M_{k}(x)}\right)^{k}\right].
\end{gather}
As a result of the asymptotic relation (\ref{eq:singlelayer_moment_asymptotic}), we expect the accuracy of $t_{s}^{(3)}(x)$ to increase as $k$ increases, and this hypothesis is tested in the Section \ref{sec:results}.

Assuming $\delta = 10^{-p}$, where $p > 0$, gives the following alternative form:
\begin{gather}
\label{eq:singlelayer_local_transition_time_3_p}
t_{s}^{(3)}(x) := \frac{\log(\alpha_{k}(x)) + p\log(10)}{\beta_{k}(x)},
\end{gather}
which leads to the conclusion that for fixed $k$ and increasing $p$, the local transition time estimate $t_{s}^{(3)}(x)$ increases linearly with slope $\log(10)/\beta_{k}(x)$, e.g., the additional transition time required when decreasing $\delta$ from $10^{-(p-1)}$ to $10^{-p}$ is equal to the additional time when decreasing $\delta$ from $10^{-p}$ to $10^{-(p+1)}$. 

Interestingly, $k = 1$ gives $\alpha_{1}(x) = 1$ and $\beta_{1}(x) = 1/M_{1}(x)$ since $M_{0}(x) = 1$. In this case, the approximation (\ref{eq:CDF_moments}) simplifies to $F(t; x) \simeq 1 - e^{-t/M_{1}(x)}$ for large $t$, which is nothing more than the cumulative distribution function of the exponential distribution with mean $M_{1}(x)$. For this simplest of cases, the estimate of the local transition time (\ref{eq:singlelayer_local_transition_time_3b}) reduces to the mean action time multiplied by a correction factor depending on the tolerance: $t_{s}^{(3)}(x) := M_{1}(x)\log(\delta^{-1}) \equiv \text{MAT}(x)\log(\delta^{-1})$.

\section{Results and Discussion}
\label{sec:results}
To investigate the accuracy of the transition time estimates presented and developed in sections \ref{sec:singlelayer_low_accuracy_estimates} and \ref{sec:singlelayer_high_accuracy_estimates}, we consider the following three test cases:

\medskip
\noindent Case A:
\begin{gather*}
D = 1.0,\quad u_{0}(x) = 0,\\
u(0,t) = 1,\quad \frac{\partial u}{\partial x}(1,t) = 0.
\end{gather*}
Case B:
\begin{gather*}
D = 0.01,\quad u_{0}(x) = 1,\\
u(0,t) - 0.1\frac{\partial u}{\partial x}(0,t) = 0,\quad u(1,t) = 0.5.
\end{gather*}
Case C:
\begin{gather*}
D = 0.1,\quad u_{0}(x) = \begin{cases} 1 & \text{if $0.25<x<0.75$},\\ 0 & \text{else},\end{cases}\\
\frac{\partial u}{\partial x}(0,t) = 0,\quad \frac{\partial u}{\partial x}(1,t) = 0.
\end{gather*} %H(x-0.25) - H(x-0.75)
where in each case $[l_{0},l_{m}] = [0,1]$. Together, the above problems test each of the three types of boundary conditions (Dirichlet, Neumann, Robin). Cases A and B resemble classical problems in heat conduction. For example, in Case A, a bar initially at temperature zero is suddenly heated at its left-boundary, how long does it take for the whole bar to heat up? On the other hand, Case C is a typical problem in mathematical biology \cite{fernando_2010}, where a region initially fully occupied by cells ($u(x,0) = 1$ for $0.25<x<0.75$) is left to diffuse. 

Each of the test cases satisfy the requirement that $F(t;x)$ (\ref{eq:singlelayer_CDF}) defines a cumulative distribution function since, at each position $x$, the solution $u(x,t)$ is either non-decreasing for all $t > 0$ if $u_{\infty}(x) > u_{0}(x)$ (as in Case A for all $0<x<1$ and Case C for $0<x<0.25$ and $0.75<x<1$) or non-increasing for all $t > 0$ if $u_{\infty}(x) < u_{0}(x)$ (as in Case B for all $0<x<1$ and Case C for $0.25<x<0.75$). For Case C, if the region centered around $x = 0$ that is fully occupied by cells, namely $0.25<x<0.75$, is shortened, however, the monotonicity property is violated.

\def\fwidth{0.25\textwidth}
\definecolor{color1}{rgb}{0,0.447,0.741}
\definecolor{color2}{rgb}{0.85,0.325,0.098}
\definecolor{color3}{rgb}{0.466,0.674,0.188}
\definecolor{color4}{rgb}{0.929,0.694,0.125}
\definecolor{color5}{rgb}{0.494,0.184,0.556}
\definecolor{color6}{rgb}{0.8,0.8,0.8}

Recall the three global transition time estimates presented in this paper, namely $\widehat{t}_{s}^{(n)}$ ($n = 1,2,3$), the definitions of which are reiterated below:
\begin{itemize}
\item $\widehat{t}_{s}^{(1)}$ [Equations \ref{eq:singlelayer_local_transition_time_1} and \ref{eq:singlelayer_global_transition_time_estimate}]\\ Based on using the mean or first moment of the probability distribution $f(t; x)$ (\ref{eq:singlelayer_PDF}) (or equivalently the mean action time) as an estimate of the local transition time, as described in section \ref{sec:singlelayer_mean_action_time}.
\item $\widehat{t}_{s}^{(2)}$ [Equations \ref{eq:singlelayer_local_transition_time_2} and \ref{eq:singlelayer_global_transition_time_estimate}]\\ Based on using the sum of the mean and standard deviation of the probability distribution $f(t; x)$ (\ref{eq:singlelayer_PDF}) (or equivalently the mean plus one standard deviation of action time) as an estimate of the local transition time, as described in section \ref{sec:singlelayer_mean_plus_standard_deviation} and previously by \citet{simpson_2013} and \citet{jazaei_2014}. 
\item $\widehat{t}_{s}^{(3)}$ [Equations \ref{eq:singlelayer_local_transition_time_3b} and \ref{eq:singlelayer_global_transition_time_estimate}]\\ Based on using the higher order moments of the probability distribution $f(t;x)$ (\ref{eq:singlelayer_PDF}) to approximate the large time asymptotic behaviour of the cumulative distribution function $F(t; x)$ (\ref{eq:singlelayer_CDF})
as described in section \ref{sec:singlelayer_high_accuracy_estimates} 
\end{itemize}

For each global transition time estimate, the maximisation problem implied by equation (\ref{eq:singlelayer_global_transition_time_estimate}) is solved by first noting that equation (\ref{eq:singlelayer_global_transition_time_estimate}) is equivalent to $\widehat{t}_{s}^{(n)} = -\min_{x\in\mathcal{R}}[-t_{s}^{(n)}(x)]$ and then using MATLAB's \texttt{fminbnd} function with option \texttt{TolX = 1e-14} \cite{matlab_fminbnd}. To calculate the moments we have implemented Algorithm \ref{alg:singlelayer_moments} in MATLAB using the Symbolic Math Toolbox \cite{matlab_symbolic}.

Recall that $\widehat{t}_{s}^{(3)}$ depends on two parameters that are free-to-choose: the prescribed tolerance $\delta$ and the index $k$ which specifies which two consecutive moments (i.e., $M_{k-1}(x)$ and $M_{k}(x)$) are utilised. Initially, we present results for $\delta = 0.02$ and $k = 2$, with the tolerance value chosen to ensure that the solution at the global transition time is visibly distinguishable, if only slightly, from the steady state $u_{\infty}(x)$. With $k = 2$, the local transition time estimate $t_{s}^{(3)}(x)$ (\ref{eq:singlelayer_local_transition_time_3b}) simplifies to:
\begin{gather*}
t_{s}^{(3)}(x) = \frac{M_{2}(x)}{2M_{1}(x)}\log\left(\frac{2M_{1}(x)^{2}}{M_{2}(x)\delta}\right),
\end{gather*}
which is comparable to $t_{s}^{(2)}(x)$ (\ref{eq:singlelayer_local_transition_time_2}) in its simplicity and dependence on the first and second moments only.

\begin{figure}
\includegraphics[width=\textwidth]{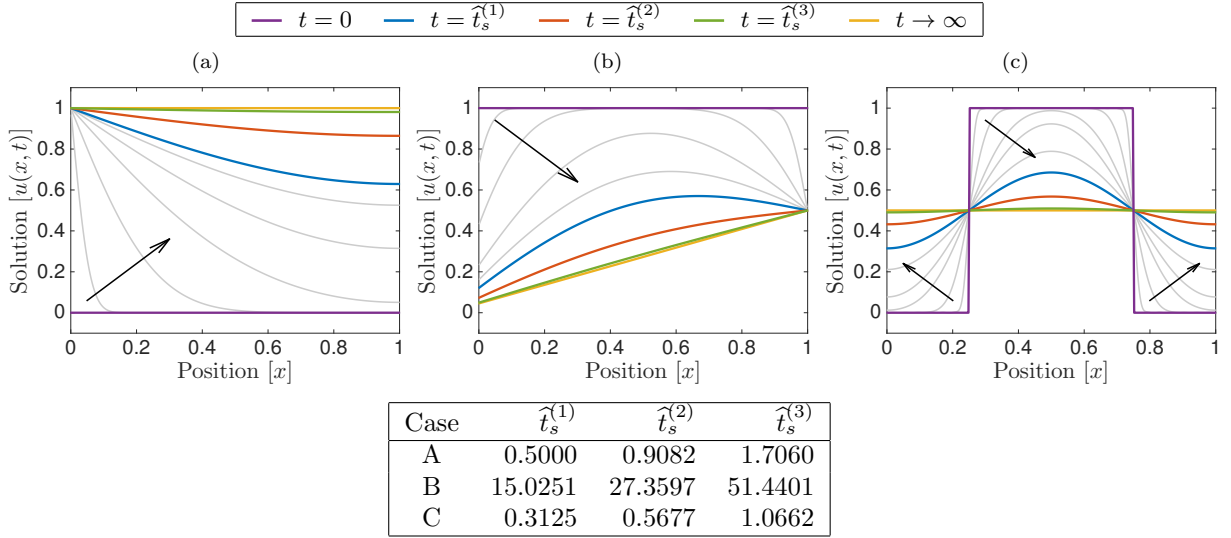}
\caption{Plot of the transient solution $u(x,t)$ of the diffusion problem (\ref{eq:singlelayer_problem}) for (a) Case A, (b) Case B and (c) Case C, depicting the transition from initial to steady state as well as the solution at each of the global transition time estimates (defined in equations \ref{eq:singlelayer_global_transition_time_estimate}, \ref{eq:singlelayer_local_transition_time_1}, \ref{eq:singlelayer_local_transition_time_2} and \ref{eq:singlelayer_local_transition_time_3b}). The value of $\widehat{t}_{s}^{(3)}$ is calculated using $\delta = 0.02$ and $k=2$. Arrows indicate the direction of increasing time.} 
\label{fig:singlelayer_solution_global_transition_time}
\end{figure}

Figure \ref{fig:singlelayer_solution_global_transition_time} plots the solution $u(x,t)$ of the diffusion problem (\ref{eq:singlelayer_problem}) for each of the three test cases, depicting the transition from initial to steady state. In these plots, $u(x,t)$ is also given at each of the global transition times estimates, with the calculated values of $t = \widehat{t}_{s}^{(n)}$ for $n = 1,2,3$ (rounded to four decimal places) tabulated in Figure \ref{fig:singlelayer_solution_global_transition_time}. In addition to the visual comparison provided by Figure \ref{fig:singlelayer_solution_global_transition_time}, in Table \ref{tab:singlelayer_global_transition_time_errors} we include errors:
\begin{gather}
\label{eq:singlelayer_residual_error}
\varepsilon_{s}^{(n)} = \max_{x\in \mathcal{R}}\left[\frac{u(x,\widehat{t}_{s}^{(n)})-u_{\infty}(x)}{u_{0}(x)-u_{\infty}(x)}\right],\quad n = 1,2,3,
\end{gather}
as a quantitive measure of how close the transient solution at each estimate is to the steady-state solution. The transient solution $u(x,\widehat{t}_{s}^{(n)})$ is evaluated by taking the first 50 terms in the classical eigenfunction expansion solution \cite{strauss_1992,trim_1990}, which is more than sufficient since the values of $\widehat{t}_{s}^{(n)}$ are relatively large. For $\widehat{t}_{s}^{(3)}$, we compute errors corresponding to three different tolerances $\delta = 0.02, 10^{-3},10^{-5}$. Noting the definition of the local transition time (Definition \ref{def:singlelayer_local_transition_time}), $\varepsilon_{s}^{(n)}$ ideally should be close to the prescribed tolerance $\delta$.

\begin{table}
\centering
\renewcommand{\tabcolsep}{3.6pt}
\begin{tabular}{|crrrrrr|}
\hline
\multirow{2}{*}{Case} & \multirow{2}{*}{$\varepsilon_{s}^{(1)}$} & \multirow{2}{*}{$\varepsilon_{s}^{(2)}$} & $\varepsilon_{s}^{(3)}$  & $\varepsilon_{s}^{(3)}$ & $\varepsilon_{s}^{(3)}$ & $\varepsilon_{s}^{(3)}$\\
& & & [$k = 2$, $\delta = 0.02$] & [$k = 2$, $\delta = 10^{-3}$] & [$k = 2$, $\delta = 10^{-5}$] & [$k = 5$, $\delta = 0.02$]\\
\hline
A & 0.3708 & 0.1354 & 0.0189 & 8.69e-04 & 7.64e-06 & 0.0200\\
B & 0.3721 & 0.1356 & 0.0188 & 8.58e-04 & 7.43e-06 & 0.0200\\
C & 0.3708 & 0.1354 & 0.0189 & 8.69e-04 & 7.64e-06 & 0.0200\\
\hline
\end{tabular}
\caption{Errors (\ref{eq:singlelayer_residual_error}), corresponding to the three different global transition time estimates (defined in equations \ref{eq:singlelayer_global_transition_time_estimate}, \ref{eq:singlelayer_local_transition_time_1}, \ref{eq:singlelayer_local_transition_time_2} and \ref{eq:singlelayer_local_transition_time_3b}), for each of the three test cases. The value of $\varepsilon_{s}^{(3)}$ is calculated using different combinations of tolerance $\delta$ and moment index $k$.}
\label{tab:singlelayer_global_transition_time_errors}
\end{table}

The following observations can be drawn from Figure \ref{fig:singlelayer_solution_global_transition_time} and Table \ref{tab:singlelayer_global_transition_time_errors}:
\begin{enumerate}
\item $\widehat{t}_{s}^{(1)}$ underestimates the time required to effectively reach steady state for all three test cases. The transition from initial to steady state is far from complete and this is confirmed by the large values of $\varepsilon_{s}^{(1)}$ in Table \ref{tab:singlelayer_global_transition_time_errors}.
\item $\widehat{t}_{s}^{(2)}$ significantly improves on $\widehat{t}_{s}^{(1)}$ (as has been reported previously by \citet{simpson_2013} and \citet{jazaei_2014} for a groundwater modelling problem), however, in all three test cases it is clearly visible that the transition from initial to steady state is still not complete.
\item $\widehat{t}_{s}^{(3)}$ uses the the same moments as $\widehat{t}_{s}^{(2)}$ but produces a far superior estimate of the global transition time with the transient solutions at $t = \widehat{t}_{s}^{(3)}$ (see Figure \ref{fig:singlelayer_solution_global_transition_time}) very close to steady state. For all three test cases, the accuracy of $\widehat{t}_{s}^{(3)}$  is quite remarkable. Using only the first and second moments, $\widehat{t}_{s}^{(3)}$ leads to errors $\varepsilon_{s}^{(3)}$ that are less than, and very close to, the prescribed tolerances of $\delta = 0.02, 10^{-3}, 10^{-5}$ (Table \ref{tab:singlelayer_global_transition_time_errors}). Taking higher consecutive moments, that is increasing $k$, further improves accuracy since the asymptotic relation (\ref{eq:singlelayer_moment_asymptotic}) is more accurate for larger $k$. This is demonstrated by choosing $k = 5$, which gives $\varepsilon_{s}^{(3)} = 0.0200 = \delta$ (to the four decimal places displayed) for all three test cases.
\end{enumerate} 

\begin{figure}[H]
\includegraphics[width=\textwidth]{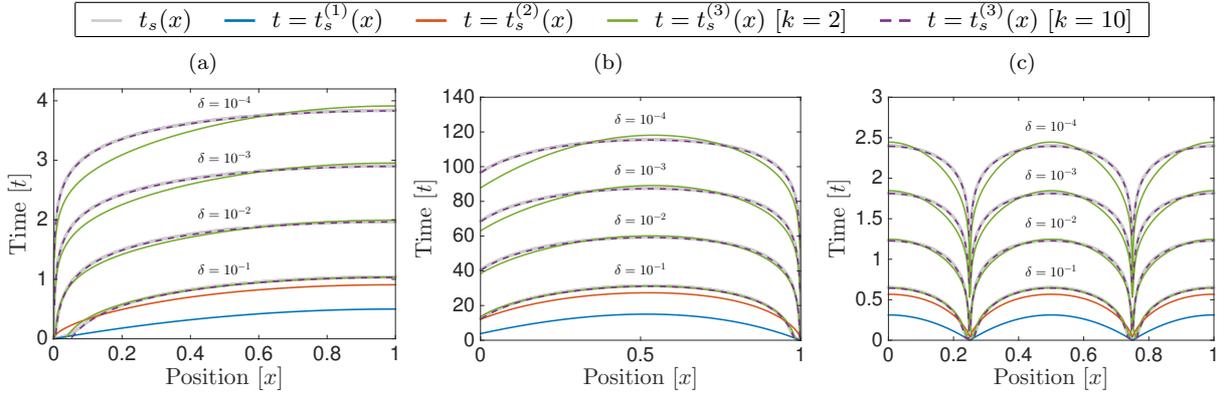}
\caption{Comparison of the exact transition time $t_{s}(x)$ (defined in equation \ref{eq:singlelayer_local_transition_time}) and local transition time estimates, $t_{s}^{(n)}(x)$ for $n = 1,2,3$ (defined in equations \ref{eq:singlelayer_global_transition_time_estimate}, \ref{eq:singlelayer_local_transition_time_1}, \ref{eq:singlelayer_local_transition_time_2} and \ref{eq:singlelayer_local_transition_time_3b}) for (a) Case A, (b) Case B and (c) Case C. The function $t_{s}^{(3)}(x)$ is shown for each combination of tolerance $\delta = 10^{-1}, 10^{-2}, 10^{-3}, 10^{-4}$ and moment index $k = 2,10$.} 
\label{fig:singlelayer_local_transition_time}
\end{figure}

In Figure \ref{fig:singlelayer_local_transition_time}, we plot the exact local transition time $t_{s}(x)$ (Definition \ref{def:singlelayer_local_transition_time}) and local transition time estimate $t_{s}^{(3)}(x)$ (\ref{eq:singlelayer_local_transition_time_3b}) for all three test cases and four different choices of the tolerance $\delta$. To calculate $t_{s}(x)$ we solve equation (\ref{eq:singlelayer_local_transition_time}) (with the first 50 terms used in the eigenfunction solution expansion for $u(x,t)$) using MATLAB's \texttt{fzero} function with the default tolerance \cite{matlab_fzero}. Recall that $t_{s}(x)$ provides the finite time required for the solution at position $x$ to transition to within a specified tolerance $\delta$ of the steady state solution at position $x$ (as defined in equation \ref{eq:singlelayer_local_transition_time}) and $t_{s}^{(3)}(x)$ is an estimate of that value. The local transition time profiles are therefore useful for determining the position(s) that take the shortest/longest time to reach steady state. For example, for Case C, the steady state is reached in the longest time at $x = 0, 0.5, 1$ and the shortest time near $x = 0.25$ and $x = 0.75$, where the solution $u(x,t) = 0.5 = u_{\infty}(x)$ for all $t > 0$. Both $t_{s}^{(1)}(x)$ and $t_{s}^{(2)}(x)$ are also useful in this regard as their shape follows closely that of $t_{s}(x)$: for both Case A and Case C, the value of $x$ that maximises/minimises $t_{s}^{(1)}(x)$ and $t_{s}^{(2)}(x)$ also maximises/minimises the exact local transition time $t_{s}(x)$. As observed previously, increasing $k$ leads to a better match with the exact transition time, with the curves almost indistinguishable for $k = 10$. Finally, one must take care when using $t_{s}^{(3)}(x)$ for moderately large values of the tolerance $\delta$. As previously remarked in section \ref{sec:singlelayer_high_accuracy_estimates}, $t_{s}^{(3)}(x)$ is non-physical if $\alpha_{k}(x) < \delta$ and this behaviour is observed in Figure \ref{fig:singlelayer_local_transition_time} for $\delta = 10^{-1}$ near $x = 0$ in Case A and near $x = 1$ in Case B. This anomaly is tied to the inaccuracy of taking only the leading term in the transient solution $u(x,t)$ (which is effectively being approximated in equation (\ref{eq:CDF_large_t_alpha_beta})) for moderately small values of time $t$. 

\begin{table}[H]
\centering
\def\arraystretch{1.2}
\renewcommand{\tabcolsep}{3.6pt}
\begin{tabular}{|lll|rrrrrr|}
\cline{4-9}
\multicolumn{3}{c|}{} & $\delta = 10^{-1}$ & $\delta = 10^{-2}$ & $\delta = 10^{-3}$ & $\delta = 10^{-4}$ & $\delta = 10^{-5}$ & $\delta = 10^{-6}$\\
\hline
 & $\widehat{t}_{s}$ & & 1.0311 & 1.9643 & 2.8975 & 3.8307 & 4.7639 & 5.6971\\
 & $\widehat{t}_{s}^{(3)}$ & [$k = 1$] & 1.1513 & 2.3026 & 3.4539 & 4.6052 & 5.7565 & 6.9078\\
 & $\widehat{t}_{s}^{(3)}$ & [$k = 2$] & 1.0354 & 1.9948 & 2.9542 & 3.9136 & 4.8730 & 5.8324\\
 & $\widehat{t}_{s}^{(3)}$ & [$k = 5$] & 1.0311 & 1.9643 & 2.8975 & 3.8308 & 4.7640 & 5.6973\\
Case A & $\widehat{t}_{s}^{(3)}$ & [$k = 10$] & 1.0311 & 1.9643 & 2.8975 & 3.8307 & 4.7639 & 5.6971\\
& $\bigl|\widehat{t}_{s} - \widehat{t}_{s}^{(3)}\bigr|/\bigl|\widehat{t}_{s}\bigr|$ & [$k = 1$] & 1.17e-01 & 1.72e-01 & 1.92e-01 & 2.02e-01 & 2.08e-01 & 2.12e-01\\
& $\bigl|\widehat{t}_{s} - \widehat{t}_{s}^{(3)}\bigr|/\bigl|\widehat{t}_{s}\bigr|$ & [$k = 2$] & 4.14e-03 & 1.55e-02 & 1.96e-02 & 2.16e-02 & 2.29e-02 & 2.38e-02\\
& $\bigl|\widehat{t}_{s} - \widehat{t}_{s}^{(3)}\bigr|/\bigl|\widehat{t}_{s}\bigr|$ & [$k = 5$] & 4.54e-05 & 2.63e-06 & 1.23e-05 & 2.05e-05 & 2.52e-05 & 2.84e-05\\
& $\bigl|\widehat{t}_{s} - \widehat{t}_{s}^{(3)}\bigr|/\bigl|\widehat{t}_{s}\bigr|$ & [$k = 10$] & 2.28e-09 & 5.79e-08 & 3.98e-07 & 5.45e-11 & 1.06e-10 & 2.39e-08\\
\hline
 & $\widehat{t}_{s}$ & & 31.0746 & 59.1707 & 87.2666 & 115.3624 & 143.4582 & 171.5541\\
 & $\widehat{t}_{s}^{(3)}$ & [$k = 1$] & 34.5967 & 69.1934 & 103.7901 & 138.3867 & 172.9834 & 207.5801\\
 & $\widehat{t}_{s}^{(3)}$ & [$k = 2$] & 31.1946 & 60.1603 & 89.1312 & 118.1046 & 147.0794 & 176.0552\\
 & $\widehat{t}_{s}^{(3)}$ & [$k = 5$] & 31.0689 & 59.1697 & 87.2706 & 115.3715 & 143.4724 & 171.5733\\
Case B & $\widehat{t}_{s}^{(3)}$ & [$k = 10$] & 31.0749 & 59.1707 & 87.2665 & 115.3624 & 143.4582 & 171.5541\\
& $\bigl|\widehat{t}_{s} - \widehat{t}_{s}^{(3)}\bigr|/\bigl|\widehat{t}_{s}\bigr|$ & [$k = 1$] & 1.13e-01 & 1.69e-01 & 1.89e-01 & 2.00e-01 & 2.06e-01 & 2.10e-01\\
& $\bigl|\widehat{t}_{s} - \widehat{t}_{s}^{(3)}\bigr|/\bigl|\widehat{t}_{s}\bigr|$ & [$k = 2$] & 3.86e-03 & 1.67e-02 & 2.14e-02 & 2.38e-02 & 2.52e-02 & 2.62e-02\\
& $\bigl|\widehat{t}_{s} - \widehat{t}_{s}^{(3)}\bigr|/\bigl|\widehat{t}_{s}\bigr|$ & [$k = 5$] & 1.83e-04 & 1.63e-05 & 4.64e-05 & 7.87e-05 & 9.86e-05 & 1.12e-04\\
& $\bigl|\widehat{t}_{s} - \widehat{t}_{s}^{(3)}\bigr|/\bigl|\widehat{t}_{s}\bigr|$ & [$k = 10$] & 8.09e-06 & 1.55e-07 & 6.07e-08 & 1.26e-08 & 1.66e-08 & 3.63e-08\\
\hline
 & $\widehat{t}_{s}$ & & 0.6444 & 1.2277 & 1.8109 & 2.3942 & 2.9774 & 3.5607\\
 & $\widehat{t}_{s}^{(3)}$ & [$k = 1$] & 0.7196 & 1.4391 & 2.1587 & 2.8782 & 3.5978 & 4.3173\\
 & $\widehat{t}_{s}^{(3)}$ & [$k = 2$] & 0.6471 & 1.2467 & 1.8464 & 2.4460 & 3.0456 & 3.6453\\
 & $\widehat{t}_{s}^{(3)}$ & [$k = 5$] & 0.6444 & 1.2277 & 1.8110 & 2.3942 & 2.9775 & 3.5608\\
Case C & $\widehat{t}_{s}^{(3)}$ & [$k = 10$] & 0.6444 & 1.2277 & 1.8109 & 2.3942 & 2.9774 & 3.5607\\
& $\bigl|\widehat{t}_{s} - \widehat{t}_{s}^{(3)}\bigr|/\bigl|\widehat{t}_{s}\bigr|$ & [$k = 1$] & 1.17e-01 & 1.72e-01 & 1.92e-01 & 2.02e-01 & 2.08e-01 & 2.12e-01\\
& $\bigl|\widehat{t}_{s} - \widehat{t}_{s}^{(3)}\bigr|/\bigl|\widehat{t}_{s}\bigr|$ & [$k = 2$] & 4.14e-03 & 1.55e-02 & 1.96e-02 & 2.16e-02 & 2.29e-02 & 2.38e-02\\
& $\bigl|\widehat{t}_{s} - \widehat{t}_{s}^{(3)}\bigr|/\bigl|\widehat{t}_{s}\bigr|$ & [$k = 5$] & 4.56e-05 & 2.78e-06 & 1.26e-05 & 2.08e-05 & 2.51e-05 & 2.83e-05\\
& $\bigl|\widehat{t}_{s} - \widehat{t}_{s}^{(3)}\bigr|/\bigl|\widehat{t}_{s}\bigr|$ & [$k = 10$] & 2.76e-07 & 2.02e-07 & 9.79e-08 & 3.20e-07 & 5.92e-08 & 7.36e-08\\
\hline
\end{tabular}
\caption{Comparison of the exact global transition time $\widehat{t}_{s}$ (defined in equation \ref{eq:singlelayer_global_transition_time}) and the global transition time estimate $\widehat{t}_{s}^{(3)}$ (defined by equations \ref{eq:singlelayer_local_transition_time_3b} and \ref{eq:singlelayer_global_transition_time_estimate}) for each combinations of tolerance $\delta = 10^{-1},10^{-2},10^{-3},10^{-4},10^{-5},10^{-6}$ and moment index $k = 1,2,5,10$.}
\label{tab:singlelayer_global_transition_time}
\end{table}

In Table \ref{tab:singlelayer_global_transition_time}, the new global transition time estimate $\widehat{t}_{s}^{(3)}$ is compared to the exact global transition time $\widehat{t}_{s}$ for several values of the tolerance $\delta$ and index $k$. Using only the moments and without explicit calculation of the transient solution $u(x,t)$, it is evident that $\widehat{t}_{s}^{(3)}$ is able to very accurately estimate $\widehat{t}_{s}$. Rough rule-of-thumb values are obtained for $k=1$, while for $k = 2$, the relative errors are all between $10^{-1}$ and $10^{-3}$ indicating an accuracy of at least one and at most three significant digits. Increasing $k$ increases the accuracy (due to the accuracy of the asymptotic relation (\ref{eq:singlelayer_moment_asymptotic}) improving) with at least 3 and and at most 6 significant digits obtained for $k = 5$, which is probably sufficient for most applications. For $k = 10$, one obtains the exact value to all four decimal places displayed with the relative errors indicating an accuracy of between 6 and 11 significant digits (inclusive).

To conclude this section, we provide some observations regarding the following generalised version of Case A, where the constant diffusivity $D$ is arbitrary and the boundary conditions are specified as:
\begin{align*}
u(0,t) = c_{L},\quad \frac{\partial u}{\partial x}(L,t) = 0.
\end{align*}
For this problem, computing the three local transition time estimates and evaluating them at $x = L$, where the maximum occurs, yields the following global transition time estimates:
\begin{align}
\label{eq:CaseD_global_transition_time_estimate_1}
\widehat{t}_{s}^{(1)} &= \frac{L^{2}}{2D},\\
\label{eq:CaseD_global_transition_time_estimate_2}
\widehat{t}_{s}^{(2)} &= \frac{L^{2}}{2D}\left(1+\frac{\sqrt{6}}{3}\right),\\
\label{eq:CaseD_global_transition_time_estimate_3}
\widehat{t}_{s}^{(3)} &= \gamma_{k}\frac{L^{2}}{D}\log\left(\frac{\theta_{k}}{\delta}\right),
\end{align}
where $\gamma_{k}$ and $\theta_{k}$ are constants that depend on the chosen moment index $k$.

Note that all three estimates are proportional to the diffusive timescale $L^{2}/D$, differing only by a multiplicative factor, which for $\widehat{t}_{s}^{(3)}$ depends on the tolerance $\delta$. As a result, for two different diffusion processes, the question of which takes longer can be answered by choosing any of the global transition time estimates (\ref{eq:CaseD_global_transition_time_estimate_1})--(\ref{eq:CaseD_global_transition_time_estimate_3}) and comparing its value for both processes. 

\begin{table}[H]
\centering
\begin{tabular}{|c|rrrr|}
\hline
$k$ & $\gamma_{k}$ & $\theta_{k}$ & $|\gamma_{k}-4/\pi^{2}|$ & $|\theta_{k}-4/\pi|$\\
\hline
2 & 0.4167 & 1.2000 & 1.14e-02 & 7.32e-02\\
4 & 0.4054 & 1.2712 & 1.60e-04 & 2.08e-03\\
6 & 0.4053 & 1.2732 & 2.03e-06 & 3.90e-05\\
8 & 0.4053 & 1.2732 & 2.51e-08 & 6.41e-07\\
10 & 0.4053 & 1.2732 & 3.10e-10 & 9.86e-09\\
12 & 0.4053 & 1.2732 & 3.83e-12 & 1.46e-10\\
14 & 0.4053 & 1.2732 & 4.72e-14 & 2.10e-12\\
16 & 0.4053 & 1.2732 & 5.83e-16 & 2.95e-14\\
18 & 0.4053 & 1.2732 & 7.20e-18 & 4.10e-16\\
20 & 0.4053 & 1.2732 & 8.89e-20 & 5.62e-18\\
\hline
\end{tabular}
\caption{Computed values of the constants $\gamma_{k}$ and $\theta_{k}$ that appear in the global transition time estimate (\ref{eq:CaseD_global_transition_time_estimate_3}) for the generalised version of Case A.}
\label{tab:theta_gamma}
\end{table}

Each of the global transition time estimates (\ref{eq:CaseD_global_transition_time_estimate_1})--(\ref{eq:CaseD_global_transition_time_estimate_3}) were calculated in Maple\footnote{with the environment variable \texttt{Digits} set to 50 \cite{maple_digits}.}, where exact fractional expressions for the constants $\gamma_{k}$ and $\theta_{k}$ appearing in equation (\ref{eq:CaseD_global_transition_time_estimate_3}) can be obtained. In Table \ref{tab:theta_gamma}, we give the corresponding values in decimal form and rounded to four decimal places to improve readability. For increasing values of $k$, observe in Table \ref{tab:theta_gamma} that $\gamma_{k}$ and $\theta_{k}$ are approaching the values $4/\pi^{2}$ ($\approx 0.4053$) and $4/\pi$ ($\approx 1.2732$). Hence, $\widehat{t}^{(3)}$ is approaching:
\begin{align*}
\frac{4}{\pi^{2}}\frac{L^{2}}{D}\log\left(\frac{4}{\pi\delta}\right),
\end{align*}
which is precisely the global transition time derived by using the leading term of the transient solution, namely
\begin{align*}
u(x,t) \sim c_{L}\left[1 - \frac{4}{\pi}\exp\left(-D\frac{\pi^{2}}{4L^{2}}t\right)\sin\left(\frac{\pi}{2L}x\right)\right],
\end{align*}
in the definition of the local transition time (Definition \ref{def:singlelayer_local_transition_time}), solving equation (\ref{eq:singlelayer_local_transition_time}) for $t$ and evaluating the result at $x = L$, where the maximum occurs.

\section{Conclusions}
In summary, we have derived a simple formula, denoted by $t_{s}^{(3)}(x)$ and given in equation (\ref{eq:singlelayer_local_transition_time_3b}), for calculating a finite measure of the time required for a diffusion process to reach steady state. This formula estimates the \textit{local transition time} defined as the time required for the transient solution to transition to within a prescribed tolerance $0<\delta\ll 1$ of its steady state, at position $x$. A finite measure of the time required for the entire diffusion process to reach steady state is then obtained by evaluating $t_{s}^{(3)}(x)$ at the value of $x$ that produces a maximum. This novel formula is attractive as it (i) avoids explicit calculation of the transient solution (ii) depends only on the prescribed tolerance $\delta$ and the $(k-1)$th and $k$th moments and (iii) can be used to calculate the exact transition time to effectively any number of significant digits by increasing $k$. Our results confirm that even for $k = 2$, for which $t_{s}^{(3)}(x)$ utilises only the first and second moments (or, equivalently, the mean and variance of action time), the accuracy is quite remarkable. In all cases, the new approach comprehensively outperforms existing strategies based on the mean action time and variance of action time, which while useful at characterising the associated time-scale, significantly underestimate transition times for diffusion processes.\\
\indent This paper has presented a proof-of-concept for the simple problem of homogeneous diffusion. Future work will focus on extending the new method to heterogeneous and higher-dimensional problems as well as other transport processes (e.g. advection-diffusion).

\section*{Acknowledgments}
This research was funded by the Australian Research Council (ARC) via the Discovery Early Career Researcher Award DE150101137.
%\end{acknowledgments}

%\nocite{}

%\section*{References}
\bibliographystyle{plainnat}
\bibliography{references}

\end{document}